\documentclass[11pt]{amsart}

\newtheorem{theorem}{Theorem}[section]
\newtheorem{proposition}[theorem]{Proposition}
\newtheorem{lemma}[theorem]{Lemma}
\newtheorem{corollary}[theorem]{Corollary}

\theoremstyle{definition}

\theoremstyle{remark}

\numberwithin{equation}{section}

\DeclareMathOperator{\tr}{Tr}

\DeclareMathOperator{\vol}{vol}

\newcommand{\inv}{^{-1}}
\newcommand{\del}{\partial}
\newcommand{\dbar}{\overline{\partial}}
\newcommand{\pair}[1]{\langle #1\rangle}

\newcommand{\bbZ}{{\mathbb Z}}
\newcommand{\bbR}{{\mathbb R}}

\newcommand{\calH}{{\mathcal H}}

\newcommand{\calL}{\mathcal L}

\newcommand{\cinf}{C^\infty}

\newcommand{\ord}[1]{^{(#1)}}
\newcommand{\norm}[1]{\left\Vert #1 \right\Vert}

\newcommand{\dtheta}{\partial_\theta}
\newcommand{\Lk}{L^{\otimes k}}
\newcommand{\lambdak}{\lambda^{(k)}}
\renewcommand{\Box}{{\mathcal D}}

\begin{document}

\title[Spectral Density Function]{The spectral density function for the Laplacian on high tensor powers of a line
bundle.}
\author{David Borthwick}
\address{Department of Mathematics and Computer Science\\
Emory University\\Atlanta}
\email{davidb@mathcs.emory.edu}
\thanks{First author supported in part by an NSF postdoctoral fellowship.}
\author{Alejandro Uribe}
\address{Mathematics Department\\
University of Michigan\\Ann Arbor}
\email{uribe@math.lsa.umich.edu}
\thanks{Second author supported in part by NSF grant DMS-0070690.}
\date{January, 2001}

\maketitle

\tableofcontents

\section{Introduction}

Let $X$ be a compact $2n$-dimensional almost K\"ahler manifold, with symplectic form $\omega$ and
almost complex structure $J$.  \textit{Almost K\"ahler} means that $\omega$ and
$J$ are compatible  in the sense that 
$$
\omega(Ju,Jv) = \omega(u,v)\text{ and } \omega(\cdot,\>J\cdot) \gg 0.
$$
The combination thus defines an associated Riemannian metric $\beta(\cdot,\cdot) =
\omega(\cdot,\>J\cdot)$.  Any symplectic manifold possesses such a structure.
We will assume further that $\omega$ is `integral' in the cohomological 
sense.   This means we can find a complex hermitian line bundle $L\to X$ with hermitian connection 
$\nabla$ whose curvature is $-i\omega$.  

Recently, beginning with Donaldson's seminal paper, \cite{Don}, the 
notion of ``nearly holomorphic'' or ``asymptotically holomorphic'' sections of $\Lk$
has attracted a fair amount of attention.  Let us recall that one natural way to 
define spaces of such sections is by means of an analogue of the 
$\dbar$-Laplacian \cite{BU96, BU00}.

The hermitian structure and connection on $L$ induce corresponding
structures on $\Lk$.  In combination with $\beta$ this defines a Laplace operator $\Delta_k$
acting on $\cinf(X; \Lk)$.  Then the sequence of operators
\[
\Box_{k} = \Delta_{k}-nk
\]
has the same principal and subprincipal symbols as the d-bar Laplacian 
in the integrable case; in fact in the K\"ahler case $\Box_{k}$ {\em 
is} the $\dbar$-Laplacian.   (By K\"ahler case we mean not only that $J$ is integrable 
but also that $L$ is hermitian \textit{holomorphic} with $\nabla$ the induced connection.)
The large $k$ behavior of the spectrum of $\Delta_k$ was 
studied (in somewhat greater generality) by  Guillemin-Uribe \cite{GU}.  For our purposes,
the main results can be summarized as follows:
\begin{theorem}\label{GUthm}\cite{GU}
There exist constants $a>0$ and $M$ (independent of $k$), such that for large $k$ the spectrum of
$\Box_{k}$ lies in $(ak, \infty)$ except for a finite number of eigenvalues contained
in $(-M,M)$.  The number $n_k$ of eigenvalues in $(-M,M)$ is a polynomial in $k$ with asymptotic
behavior $n_k \sim k^n \vol(X)$, 
and can be computed exactly by a symplectic Riemann-Roch formula.

Furthermore, if the eigenvalues in $(-M,M)$ are labeled $\lambdak_j$, then there exists a
\textit{spectral density function} $q\in \cinf(X)$ such that for any $f\in C(\bbR)$,
$$
\frac1{n_k} \sum_{j=1}^{n_k} f(\lambdak_j) \longrightarrow \frac1{\vol(X)}
\int_X (f\circ q) \frac{\omega^n}{n!},
$$ 
as $k\to \infty$.
\end{theorem}
The proof of Theorem \ref{GUthm} is based on the analysis of generalized Toeplitz structures
developed in \cite{BG}.

By the remarks above, in the K\"ahler case all 
$\lambdak_j = 0$, corresponding to eigenfunctions which are holomorphic sections of $\Lk$.  
Hence $q\equiv 0$ for a true K\"ahler structure.  In general, it is 
therefore natural to consider sections of $\Lk$ spanned by the 
eigenvalues of $\Box_{k}$ in $(-M,M)$ as being analogous to 
holomorphic sections.

The goal of the present paper is to derive a simple geometric formula for the spectral
density function $q$. Our main result is:
\begin{theorem}\label{mainthm}
The spectral density function is given by
$$
q =  - \frac5{24} |\nabla J|^2
$$
\end{theorem}

\begin{corollary}
The spectral density function is identically zero iff $(X,J, \omega)$ 
is K\"ahler.
\end{corollary}

It is natural to ask if one can choose $J$ so that $q$ is very small, 
i.\,e.\ if the symplectic invariant
\[
j(X,\omega) := \inf\{\,\norm{|\nabla J|^2}_{\infty}\;;\; J\ \text{a compatible 
almost complex structure}\ \}
\]
is always zero.  We have learned from Miguel Abreu that for 
Thurston's manifold $j=0$; it would be very interesting to find 
$(X,\omega)$ with $j>0$.

The proof of Theorem \ref{mainthm} starts with the standard and very useful
observation that sections of $\Lk$ are equivalent to equivariant functions on
an associated principle bundle $\pi:Z\to X$.  We endow $Z$ with a `Kaluza-Klein' metric
such that the fibers are geodesic.  Then the main idea exploited in the proof is the
construction of approximate eigenfunctions (quasimodes) of the Laplacian $\Delta_Z$
concentrated on these closed geodesics.  Such quasimodes are equivariant and thus
naturally associated to sections of $\Lk$.  Moreover, the value of the spectral density function 
$q(x)$ is encoded in the eigenvalue of the quasimode concentrated on the fiber 
$\pi^{-1}(x)\subset Z$.  

\section{Preliminaries}\label{prelim}

The associated principle bundle to $L$ is easily obtained as the unit circle bundle
$Z \subset L^*$.  There is a 1-1 correspondence between sections of $\Lk$ and functions on $Z$
which are $k$-equivariant with respect to the $S^1$-action, i.e. $f(z.e^{i\theta}) =
e^{ik\theta} f(z)$.

The connection $\nabla$ on $L$ induces a connection 1-form $\alpha$ on $Z$.  The curvature
condition on $\nabla$ translates to
$$
d\alpha = \pi^* \alpha,
$$
where $\pi:Z\to X$.  
Together with the Riemannian metric on $X$ and the standard metric on $S^1 = \bbR/2\pi\bbZ$,
this defines a `Kaluza-Klein' metric $g$ on $Z$ such that the projection $Z\to X$ is a
Riemannian submersion with totally geodesic fibers.  
With these choices the correspondence between equivariant functions and sections extends to
an isomorphism between 
\begin{equation}\label{lklzk}
L^2(X,\Lk)\simeq L^2(Z)_k,
\end{equation}
where $L^2(Z)_k$ is the $k$-th isotype of $L^2(Z)$ under the $S^1$ action.  

The Laplacian on $Z$ is denoted $\Delta_Z$.  By construction it commutes with the generator
$\dtheta$ of the circle action, and so it also commutes with the `horizontal Laplacian':
\begin{equation}
\Delta_h = \Delta_Z + \dtheta^2.
\end{equation}
The action of $\Delta_h$ on $L^2(Z)_k$ is equivalent under (\ref{lklzk}) to the
action of $\Delta_k$ on $L^2(X,\Lk)$.

For sufficiently large $k$, we let $\calH_k\subset L^2(Z)_k$ denote the span of the eigenvectors
with eigenvalues in the bounded range $(-M,M)$.  The corresponding orthogonal projection is denoted
$\Pi_k:L^2(Z)\to \calH_k$.  The following fact appears in the course of the proof of Theorem \ref{GUthm}:
\begin{lemma}\label{qlemma}\cite{GU}
There is a sequence of functions $q_j\in \cinf(X)$ such that
$$
\Bigl\Vert\Pi_k \Bigl(\Delta_h - nk - \sum_{j=0}^{N} k^{-j} \pi^*q_j\Bigr) \Pi_k \Bigr\Vert 
= O(k^{-(N+1)}).
$$
Moreover, the spectral density function $q$ in Theorem \ref{GUthm}
is equal to $q_0$.
\end{lemma}

\section{Quasimodes on the circle bundle}

The key to the calculation of the spectral density function at $x_0 \in X$ is the observation
that, with the Kaluza-Klein metric, the assumptions on $X$ imply the stability of the geodesic fiber 
$\Gamma = \pi^{-1}(x_0)$.  Thus one should be able to construct an
approximate eigenfunction, or \textit{quasimode}, for $\Delta_Z$ which is asymptotically 
localized on $\Gamma$.  The lowest eigenvalue of the quasimode (or rather a particular
coefficient in its asymptotic expansion) will yield the spectral density function.

The computation is largely a matter of interpolating between two natural coordinate systems.
From the point of view of writing down the Kaluza-Klein metric explicitly, the obvious
coordinate system to use is given by first trivializing $Z$ to identify a neighborhood of
$\Gamma$ with $S^1 \times U_{x_0}$, where $U_{x_0}$ is a neighborhood of $x_0$ in $X$ (the base point 
$x_0$ will be fixed throughout this section).  On $U_{x_0}$ we can introduce geodesic normal coordinates
centered at $x_0$.  These coordinates will be denoted $(\theta, x^1,\dots, x^{2n})$.  The
base point $z_0 \in \Gamma$ corresponding to $\theta=0$ is arbitrary. 
In such coordinates the connection $\alpha$ takes the form $\alpha = d\theta + \alpha_j dx^j$.

We will follow the quasimode construction outlined in Babich-Buldyrev \cite{Ba}, which is essentially
based in the normal bundle $N\Gamma \subset TZ$.   
Let $\psi:N\Gamma\to Z$ be the map defined on each fiber $N_z\Gamma$ by the
restriction of the exponential map $\exp_z:T_zZ\to Z$.  Of course, $\psi$ is only a diffeomorphism near
$\Gamma$. The \textit{Fermi coordinate system} along $\Gamma$ is defined by the combination of $\psi$
and the choice of a parallel frame for $N\Gamma$.  Let $\gamma(s)$ be a parametrization of
$\Gamma$ by arclength, with $\gamma(0) = z_0$, $\gamma'(0) = \dtheta$.  Let $e_j(s)$ be the frame
for $N_{\gamma(s)}\Gamma$ defined by parallel transport from the initial value $e_j(0) =
\del_j$.  Then the Fermi coordinates are defined by 
$$
(s,y^j) \mapsto \psi(y^j e_j(s)).
$$
Note that $s = \theta$ only on $\Gamma$.

\subsection{The ansatz}
Now we can formulate the construction of an asymptotic eigenfunction as a set
of parabolic equations on $N\Gamma$.  
Let $\kappa$ be an asymptotic parameter (eventually to be related to
$k$).  Setting $u = e^{i\kappa s} U$ we consider the equation
\begin{equation}\label{eiksu}
(\Delta_Z - \lambda)e^{i\kappa s}U(s,y) = 0.
\end{equation}
Since we are hoping to localize near $y=0$, the ansatz is to substitute 
$u^j = \sqrt{\kappa} \>y^j$ and do a formal expansion
\begin{equation}\label{DZexp}
e^{-i\kappa s} \Delta_Z e^{i\kappa s} = \kappa^2 + \kappa \calL_0
+ \sqrt{\kappa} \calL_1 + \calL_2 + \dots.
\end{equation}
This defines differential operators $\calL_j$ on a neighborhood of the zero-section in $N\Gamma$,
but since the coefficients are polynomial in the $y^j$ variables, they extend naturally to all of
$N\Gamma$.  We also make an ansatz of formal expansions for $\lambda$ and $U$:
\begin{equation*}\begin{split}
\lambda &= \kappa^2 +\sigma + \dots\\
U &= U_0 + \kappa\inv U_1 + \dots
\end{split}\end{equation*}
Substituting these expansions into (\ref{eiksu}) and reading off the orders gives the
equations
\begin{equation}\label{lueqns}\begin{split}
\calL_0 U_0 &= 0\\
\calL_1 U_0 &= 0\\
\calL_0 U_1 &= -(\calL_2 - \sigma) U_0
\end{split}\end{equation}
Since $\calL_j$ is well-defined on $N\Gamma$, we can seek global solutions $U_j(s,y)$, subject
to the boundary condition $\lim_{|y|\to\infty} U_j = 0$.
It turns out that $\calL_0$ is a very familiar parabolic operator and the second equation is satisfied
as a trivial consequence of the first.   Solutions of the third equation exist only for a certain value of
$\sigma$, and the main goal of this section is to compute this quantity.

By pulling back by $\psi$, we can use $(\theta, x)$ as an alternate coordinate system on
$N\Gamma$ (near the zero section).   We'll use $\bar\beta_{ij}, \bar\alpha_i, \bar\omega_{ij},
\bar J^i_j$ to denote the various tensors lifted from $X$ and written in these coordinates
(so all are independent of $\theta$).  
We let $\bar\Gamma_{\mu\nu}^\sigma$ denote the Christoffel symbols of the Kaluza-Klein metric $g$ 
in the $(\theta,x)$ coordinates.
To reduce notational complexity insofar as possible, we will adopt the convention that unbarred expressions
involving $\beta_{ij},
\alpha_i, \omega_{ij}, J^i_j$ and their derivatives are to be evaluated at the base point 
$x_0\in X$, e.g. $\beta_{ij} = \bar\beta_{ij}|_{x=0}$.   
The Christoffel symbols of $\beta_{ij}$ (evaluated at $x_0$) will be denoted by $F_{jk}^l$, with the
same convention for derivatives. 
(Thus $F_{jk}^l = 0$, but its derivatives are not zero.)
The freedom in the trivialization of $Z$ may be exploited to assume that
$$
\alpha_j = 0, \qquad \del_j \alpha_k = \frac12 \omega_{jk}.
$$ 

We'll use $g_{\mu\nu}$ to denote the Kaluza-Klein metric expressed in the $(\theta,x)$ 
coordinates (with the convention that Greek indices range over $0,\dots,2n$ and Roman over 
$1,\dots,2n$).  Let $\del_j$ denote the vector field
$\frac\del{\del x^j}$ on $X$.  The horizontal lift of $\del_j$ to $Z$ is then
\begin{equation}\label{ejdef}
E_j = \del_j - \bar\alpha_j \dtheta.
\end{equation}
The Kaluza-Klein metric is determined by the conditions:
$$
g(E_j,\dtheta)=0,\qquad g(\dtheta,\dtheta)=1, \qquad g(E_j,E_k) = \bar\beta_{jk}.
$$
Substituting in with (\ref{ejdef}) we quickly see that
$$
g_{00} = 1,\qquad g_{j0} = \bar\alpha_j, \qquad g_{jk} = \bar\beta_{jk} +
\bar\alpha_j\bar\alpha_k.
$$
In block matrix form we can write
\begin{equation}\label{kkg}
g = \begin{pmatrix}1&\bar\alpha\\ \bar\alpha&\bar\beta + \bar\alpha \bar\alpha\\ \end{pmatrix},
\end{equation}
from which
\begin{equation}
g\inv = \begin{pmatrix}1 + \bar\alpha\bar\beta\inv\bar\alpha&-\bar\beta\inv\bar\alpha\\
-\bar\beta\inv\bar\alpha& \bar\beta\inv \\
\end{pmatrix}.
\end{equation}

We'll use $G_{\mu\nu}$ to denote the Kaluza-Klein metric written in the Fermi
coordinates $(s,y)$, i.e.
$$
G_{00} = g(\tfrac{\del}{\del s}, \tfrac\del{\del s}), 
\qquad G_{0j} = g(\tfrac{\del}{\del s}, \tfrac\del{\del y^j}), 
\qquad G_{ij} = g(\tfrac{\del}{\del y^i}, \tfrac\del{\del y^j})
$$
$G_{\mu\nu}$ is well-defined in a neighborhood of $y=0$, and with the ansatz above we only need
to know its Taylor series to determine $\calL_j$.  As noted, the heart of the calculation will be the
change of coordinates from $(\theta,x)$ to $(s,y)$.

By assumption $G_{\mu\nu} = \delta_{\mu\nu}$ to second order in $y$.  After the substitution
$u_j = \sqrt{\kappa} y_j$, we can
write the Taylor expansions of various components as 
\begin{equation}\label{Gtaylor}\begin{split}
G_{00} &= 1 + \kappa^{-1} a\ord2 + \kappa^{-3/2} a\ord3 + \kappa^{-2} a\ord4 + \dots\\
G_{0j} &= \kappa^{-1} b_j\ord2 + \kappa^{-3/2} b_j\ord3 + \dots\\
G_{jk} &= \delta_{jk} + \kappa^{-1} c\ord2_{jk} + \dots,\\
\end{split}\end{equation}
where $(l)$ denotes the term which is a degree $l$ polynomial in $u$.  
Then using the definition
$$
\Delta_Z = - \frac1{\sqrt{G}} \del_\mu \Bigl[\sqrt{G} G^{\mu\nu} \del_\nu \Bigr],
$$
we can substitute the expansions (\ref{Gtaylor}) into (\ref{DZexp}) and read off the first few
orders in $\kappa$:
\begin{equation}\label{calls}
\begin{split}
\calL_0 &= - 2i\del_s - a\ord2 - \del_u^2\\
\calL_1 &= - a\ord3 + 2i {b^j}\ord2 \tfrac\del{\del u^j} + i (\tfrac\del{\del u^j}{b^j}\ord3)\\
\calL_2 &= - \del_s^2 + 2ia\ord2 \del_s  - a\ord4 + (a\ord2)^2 + (b\ord2)^2 \\
&\qquad + i\Bigl[ -\frac{1}2 \del_s \tr c\ord2 + 2 {b^j}\ord3 \tfrac\del{\del u^j} 
+ (\tfrac\del{\del u^j}{b^j}\ord3)\Bigr] \\
&\qquad + {c^{jk}}\ord2 \tfrac\del{\del u^j}\tfrac\del{\del u^k} 
+ (\tfrac\del{\del u^j} {c^{jk}}\ord2)\tfrac\del{\del u^k}
- \frac12\tfrac\del{\del u^j} [a\ord2 + \tr c\ord2] \tfrac\del{\del u^j}
\end{split}
\end{equation}

\subsection{The metric in Fermi coordinates}

For use in the calculation, let us first work out some simple implications of $J^2 = -1$.
Using conventions as above, this means $\bar J_j^k \bar J_k^m = -\delta_j^m$.  Differentiating
at the base point $x_0$ gives us
$$
(\del_l J_j^k) J_k^m = - J_j^k (\del_l J_k^m), \qquad
J_j^k (\del_l J_k^j) = 0 
$$
The other basic fact is $d\omega = 0$, which translates to
$$
\del_l \omega_{jk} + \del_j \omega_{kl} + \del_k \omega_{lj} = 0.
$$
\begin{lemma}
$$
\del_l J_j^l = 0.
$$
\end{lemma}
\begin{proof}
Using the fact that $J_j^l = \omega_{jk} \beta^{kl}$ we have
\begin{equation*}\begin{split}
J_j^k (\del_l J_k^l) &= - (\del_l J_j^k) J_k^l \\
&= - (\del_l \omega_{jk}) \omega^{kl}\\
&= - \frac12 (\del_l \omega_{jk} - \del_k \omega_{jl}) \omega^{kl}\\
&= \frac12 (\del_j \omega_{kl}) \omega^{kl}\\
&= - \frac12 (\del_j J_k^l) J_l^k\\
&= 0
\end{split}\end{equation*}
\end{proof}
A similar fact, which will also be useful, is:
\begin{lemma}\label{vjlemma}
For any vector $v^j$ we have
$$
(\del_l J_j^m) v^j (\omega v)_m = 0.
$$
\end{lemma}
\begin{proof}
\begin{equation*}\begin{split}
(\del_l J_j^m) v^j (\omega v)_m &=  (\del_l J_j^m) v^j J_m^s v_s \\
&=  - (\del_l J_m^s) v^j J_j^m v_s \\
&= - (\del_l \omega_{ms}) (Jv)^m v^s\\
&= (\del_l \omega_{sm}) (Jv)^m v^s\\
&= - (\del_l J_s^m) (\omega v)_m v^s\\
\end{split}\end{equation*}
\end{proof}

To proceed, we must determine the terms in the Taylor expansion of $G_{\mu\nu}$ in terms
of the geometric data $\beta,\omega,J,\alpha$.  In terms of $\del_j$, let the parallel frame
$e_j(s)$ be written $T^k_j \del_k$ (fixing $j$ for the moment).  The parallel condition is
$$
\del_s T^k_j = - \bar\Gamma^k_{0l} T^l_j,
$$
where 
$$
\bar\Gamma^k_{0l}|_{x=0} = \frac12 \beta^{km} (\del_l\alpha_m - \del_m \alpha_l)
= \frac12 \beta^{km} \omega_{lm} = \frac12 J^k_l.
$$
The solution is
$$
T^k_j = (e^{-\frac{s}2 J})^k_j. 
$$
Since this is the matrix relating the $x$-frame to the $y$-frame at $x=0$, we have
$\frac{\del x^k}{\del y^j} = T^k_j$.  This makes it convenient to introduce
$$
z^k = T^k_j y^j.
$$

The transformation to Fermi coordinates may now be written as 
$$
\theta = s + A(s,z), \qquad x^j = z^j + B^j(s,z).
$$  
The functions $A$ and $B$ are determined by the condition that the ray $t \mapsto (s,ty)$
be a geodesic.  Of course, we are really just interested in the Taylor expansions:
\begin{equation*}\begin{split}
A &= \kappa^{-1} A\ord2 + \kappa^{-3/2} A\ord3 +  \kappa^{-2} A\ord4 +\dots,\\
B^j &= \kappa^{-1} {B^j}\ord2 + \kappa^{-3/2} {B^j}\ord3 +\dots,
\end{split}\end{equation*}
where degrees are labeled as above.

Denoting the $t$ derivative by a dot, the geodesic equations are
\begin{equation}\label{tgeo}
\begin{split}
\ddot\theta &= -\Gamma_{00}^0 \dot\theta^2 -2\Gamma_{0l}^0 \dot\theta \dot x_l 
- \Gamma_{jl}^0 \dot x_j \dot x_l\\
\ddot x_k &= -\Gamma_{00}^k \dot\theta^2 -2\Gamma_{0l}^k \dot\theta \dot x_l 
- \Gamma_{jl}^k \dot x_j \dot x_l
\end{split}
\end{equation}
The Christoffel symbols of $g_{ij}$ are
\begin{equation*}
\begin{split}
&\bar\Gamma_{00}^0 = \bar\Gamma_{00}^j =0 \\
&\bar\Gamma_{0j}^0 = \frac12 (\bar J\bar\alpha)_j \\
&\bar\Gamma_{jk}^0 = \frac12 \Bigl[ \del_j\bar\alpha_k + \del_k\bar\alpha_j 
+ \bar\alpha_j(\bar J\bar\alpha)_k + \bar\alpha_k(\bar J\bar\alpha)_j \Bigr] 
- \bar F_{jk}^l \bar\alpha_l\\
&\bar\Gamma_{0k}^j = - \frac12 \bar J_k^j\\
&\bar\Gamma_{lk}^j = - \frac12 \bar J_l^j \bar\alpha_k - \frac12 \bar J_k^j \bar\alpha_l +
\bar F_{lk}^j
\end{split}
\end{equation*}

Substituting the Taylor expansion of the Christoffel symbols into (\ref{tgeo}) 
and equating coefficients, we find $A\ord2 = 0$, $B\ord2 = 0$, 
\begin{equation}\label{aab}\begin{split}
A\ord3 &= -(\del_m \del_j \alpha_l) z^m z^j z^l\\
A\ord4 &= - \frac1{24} (\del_k \del_m \del_j\alpha_l) z^k z^m z^j z^l
- \frac1{24} (\del_k F^i_{jl}) z^k z^j z^l (\omega z)_i,\\
{B^k}\ord3 &= -\frac16 (\del_m F^k_{jl}) z^m z^j z^l.
\end{split}\end{equation}
Using $x = z + \kappa^{-3/2}B\ord3 + \dots$, we can then determine the coefficients of the
expansion of $\bar\alpha_k$:
\begin{equation}\label{alphas}\begin{split}
\bar\alpha_k\ord1 &=  -\frac12 (\omega z)_k \\
\bar\alpha_k\ord2 &=  \frac12 (\del_l \del_m \alpha_k) z^l z^m \\
\bar\alpha_k\ord3 &=  \frac16 (\del_j \del_l\del_m \alpha_k) z^j z^l z^m + \frac1{12}
\omega_{ki}  (\del_m F^i_{jl}) z^m z^j z^l
\end{split}\end{equation}

The Fermi coordinate vector fields are
\begin{equation*}
\begin{split}
\del_s &= (1+ \del_s A)\del_0 + ({z'}^l + {B'}^l)\del_l,\\
\tfrac\del{\del y^j} &= (\tfrac\del{\del y^j}A) \del_0 + (T^l_j + \tfrac\del{\del
y^j}B^l)\del_l.
\end{split}
\end{equation*}
Note that $z^j = T^j_k(s) y^k$, so ${z'}^j = - \frac12(Jz)^j$.
To compute $a\ord{l}$, we use (\ref{aab}) and (\ref{alphas}) to expand 
$G_{00} = g(\del_s,\del_s)$.  The second order term is 
\begin{equation}
a\ord2 = 2\alpha_l\ord1 {z'}^l + {z'}^l z'_l = -\frac{z^2}4
\end{equation}
At third order we have
\begin{equation*}
\begin{split}
a\ord3 &= 2(A\ord3)' + 2\alpha_m\ord2 {z'}^m\\
&= -\frac13 (\del_j \del_l \alpha_m) [2{z'}^j z^l z^m + z^j z^l {z'}^m] + (\del_j \del_l \alpha_m)
z^j z^l {z'}^m\\
&= \frac 13 (\del_j \del_l \alpha_m) (Jz)^j z^l z^m - \frac 13 (\del_j \del_l \alpha_m)
z^j z^l (Jz)^m \\
&= - \frac13 (\del_l \omega_{jm}) z^j z^l (Jz)^m
\end{split}
\end{equation*}
Thus, by Lemma \ref{vjlemma} we have
\begin{equation}
a\ord3 = 0.
\end{equation}

The fourth order term is somewhat more complicated:
$$
a\ord4 = 2{A'}\ord4 + 2\alpha_m\ord3 {z'}^m + 2\alpha_m\ord1 ({B'}^m)\ord3
+ {z'}^l (\beta_{lm}\ord2 + \alpha_l\ord1 \alpha_m\ord1) {z'}^m + 2{z'}^m ({B'}^m)\ord3
$$
We'll expand the first term,
\begin{equation*}\begin{split}
2{A'}\ord4 &=  \frac1{24} (\del_k \del_m \del_j\alpha_l) [3z^k z^m (Jz)^j z^l 
+ z^k z^m z^j (Jz)^l] \\
&\qquad + \frac1{24} (\del_k F^i_{jl}) [(Jz)^k z^j z^l (\omega z)_i + 2 z^k z^j (Jz)^l
(\omega z)_i + z^k z^j z^l z_i]
\end{split}\end{equation*}
and the second,
$$
2\alpha_k\ord3 {z'}^k = - \frac16 (\del_j \del_l\del_m \alpha_k) z^j z^l z^m (Jz)^k -
\frac1{12} \omega_{ki}  (\del_m F^i_{jl}) z^m z^j z^l (Jz)^k.
$$
The terms involving $\del_m \alpha_k$ combine to form factors of $\omega_{mk}$:
\begin{equation*}\begin{split}
2{A'}\ord4 + 2\alpha_k\ord3 {z'}^k &= -\frac18 (\del_j \del_l \omega_{mk}) z^j z^l z^m (Jz)^k
+ \frac1{24} (\del_k F^i_{jl}) (Jz)^k z^j z^l (\omega z)_i \\
&\qquad + \frac1{12} (\del_k F^i_{jl}) z^k z^j (Jz)^l (\omega z)_i
+ \frac18 (\del_k F^i_{jl}) z^k z^j z^l z_i
\end{split}\end{equation*}
After noting that $2\alpha_m\ord1 ({B'}^m)\ord3+ 2{z'}^m ({B'}^m)\ord3=0$, we are left
with the term 
$$
{z'}^l (\beta_{lm}\ord2 + \alpha_l\ord1 \alpha_m\ord1) {z'}^m = \frac18  (\del_j\del_k \beta_{lm})
(Jz)^l z^j z^k (Jz)^m + \frac{z^4}{16}
$$
So in conclusion,
\begin{equation}\begin{split}
a\ord4 &= -\frac18 (\del_j \del_l \omega_{mk}) z^j z^l z^m (Jz)^k
+ \frac1{24} (\del_k F^i_{jl}) (Jz)^k z^j z^l (\omega z)_i \\
&\qquad + \frac1{12} (\del_k F^i_{jl}) z^k z^j (Jz)^l (\omega z)_i
+ \frac18 (\del_k F^i_{jl}) z^k z^j z^l z_i \\
&\qquad+ \frac18  (\del_j\del_k \beta_{lm})
(Jz)^l z^j z^k (Jz)^m + \frac{z^4}{16}
\end{split}\end{equation}

For $b_j = g(\del_s, \del_{y^j})$ the third order term will prove irrelevant, so we compute only
\begin{equation}
\begin{split}
b_j\ord2 &= \del_{y^j} A\ord3 + \alpha_m\ord2 T^m_j \\
&= - \frac16 (\del_k \del_l \alpha_m) [2T^k_j z^l z^m + z^k z^l T^m_j]
+\frac12 (\del_k \del_l \alpha_m) z^k z^l T^m_j\\
&= -\frac13 (\del_k \del_l \alpha_m) T^k_j z^l z^m + \frac13 (\del_k \del_l \alpha_m) z^k z^l
T^m_j\\
&= \frac13 (\del_l\omega_{km}) z^k z^l T^m_j
\end{split}
\end{equation}

Finally, we have $c_{lm} = g(\del_{y^l},\del_{y^m})$.  It is convenient to insert factors of $T$:
\begin{equation}\begin{split}
T^l_j c_{lm}\ord2 T^m_k  &=  \beta_{jk}\ord2 + 
\alpha_j\ord1 \alpha_k\ord1 + (\del_{z^j}B_k\ord3)  + (\del_{y_k}B_j\ord3) \\
&= \frac12 (\del_l\del_m\beta_{jk}) z^l z^m + \frac14 (\omega z)_j (\omega z)_k
- \frac16 (\del_j F_{ilk}) z^i z^l\\
&\quad - \frac13 (\del_m F_{jlk}) z^m z^l
- \frac16 (\del_k F_{ilj}) z^i z^l - \frac13 (\del_m F_{klj}) z^m z^l
\end{split}\end{equation}

\subsection{Parabolic equations}
The first of the equations (\ref{lueqns}) involves the operator
$$
\calL_0 = - 2i\del_s + \frac{u^2}4 - \del_u^2
$$
The equation $\calL_0 U_0 = 0$ is then instantly recognizable as the Schr\"odinger equation for a harmonic
oscillator.  The ``ground state'' solution 
\begin{equation}\label{u0def}
U_0 =  e^{-ins/2} e^{-u^2/4}.
\end{equation}
Now $e^{i\kappa s} U$ is supposed to be periodic, which means we must require
$$
\kappa - \frac{n}2 = k\in \bbZ.
$$
A function on $z$ which is $e^{iks}\times$(periodic) comes from a section of $L^k$, so this
$k$ is our usual asymptotic parameter, and
\begin{equation}
\kappa^2 = k^2 + nk + \frac{n^2}4
\end{equation}

By the standard analysis of the quantum harmonic oscillator, a complete set of solutions to $\calL_0 U = 0$
can be generated by application of the ``creation operator'' 
$$
\Lambda^*_j = -i e^{-is/2} (\del_{u^j} - \frac{u_j}2)
$$ 
We will need
$$
U_{ij} = \Lambda^*_i \Lambda^*_j U_0, \qquad U_{ijkl} =  \Lambda^*_i \Lambda^*_j \Lambda^*_k \Lambda^*_l U_0,
$$
which are easily computed explicitly:
\begin{equation*}\begin{split}
U_{ij} &= (-u_ju_k + \delta_{ij})e^{-is}U_0,\\
U_{ijkl} &= \Bigl(u_iu_ju_ku_l - \delta_{ij}u_ku_l - \delta_{ik}u_ju_l - \delta_{il}u_ju_k
 - \delta_{jk}u_iu_l \\
&\qquad- \delta_{kl}u_iu_j - \delta_{lj}u_iu_k + \delta_{ij}\delta_{kl} + \delta_{il}\delta_{jk} +
\delta_{ik}\delta_{jl}\Bigr)e^{-2is}U_0 
\end{split}\end{equation*}

Since $a\ord3=0$ and $\del_{u^j}{b^j}\ord2= 0$, the next operator is
$$
\calL_1 = 2i {b^j}\ord2 \tfrac\del{\del u^j}
$$
It then follows from $b_j\ord2 u^j = 0$ that
$$
\calL_1 U_0 = 0.
$$
Moreover, it is easy to check, using the creation operators, that $U_0$ is the unique solution of 
$\calL_0 U = 0$ for which this is true.

Consider finally the third equation
\begin{equation}\label{l2eq}
\calL_0 U_1 = - (\calL_2 - \sigma) U_0,
\end{equation}
from which we'll determine $\sigma$.
Since $\calL_2 U_0$ has coefficients polynomial in $u_j$ of order no more than four, we can expand
\begin{equation}\label{l2u}
\calL_2 U_0 = [C^{ijkl} u_iu_ju_ku_l + C^{ij} u_iu_j+C]U_0.
\end{equation}
\begin{proposition}
The equations (\ref{lueqns}) have a solution $U_0$, $U_1 \in \cinf(N\Gamma)$ 
if and only if 
\begin{equation}\label{sigcont}
\sigma =  C + {C_l}^l + 3{C_{kk}}^{ll},
\end{equation}
where the coefficients $C^{ijkl}$ are assumed symmetrized.
\end{proposition}
\begin{proof}
In terms of the basis for the kernel of $\calL_0$ we can rewrite (\ref{l2u}) as
$$
\calL_2 U_0 = e^{2is} D^{ijkl} U_{ijkl}  + e^{is}  D^{ij} U_{ij} + D U_0.
$$ 
Observe that
$$
\calL_0 (e^{2is}U_{ijkl}) = -4U_{ijkl},\qquad \calL_0 (e^{is}U_{ij}) = -2U_{ij}.
$$
So the equation $\calL_0 U_1 = - (\calL_2-\sigma)U_0$ has a solution only if $\sigma = D$, in which case we
can set
$$
U_1 = \frac14 e^{2is} D^{ijkl} U_{ijkl} + \frac12 e^{is}  D^{ij} U_{ij}. 
$$

To compute $D$ we note
$$
C^{ij} u_iu_j U_0 = - C^{ij} e^{is} U_{ij} + {C_l}^l U_0,
$$
and (with the symmetry assumption),
\begin{equation*}\begin{split}
C^{ijkl} u_iu_ju_ku_l U_0 &= C^{ijkl} e^{2is} U_{ijkl} + \Bigl[6{C_j}^{jkl}u_ku_l 
- 3{C_{kk}}^{ll} \Bigr] U_0\\
&=C^{ijkl} e^{2is} U_{ijkl} + (\dots)e^{is} U_{jk} + 3{C_{kk}}^{ll} U_0
\end{split}\end{equation*}
This means that
$$
D =  C + {C_l}^l + 3{C_{kk}}^{ll}.
$$
\end{proof}

To conclude the computation, we will examine $\calL_2 U_0$ piece by piece and form the 
contractions of coefficients according to (\ref{sigcont}).  From 
(\ref{calls}) we break up $\calL_2 U_0 = W_1 + \dots W_6$, where
\begin{equation*}\begin{split}
W_1 &= [- \del_s^2 + 2ia\ord2 \del_s]U_0\\
W_2 &= [- a\ord4 + (a\ord2)^2]U_0\\
W_3 &= (b\ord2)^2 U_0 \\
W_4 &= i\Bigl[ -\frac{1}2 \del_s \tr c\ord2 + 2 {b^j}\ord3 \tfrac\del{\del u^j} 
+ (\tfrac\del{\del u^j}{b^j}\ord3)\Bigr]U_0 \\
W_5 &= \Bigl[{c^{jk}}\ord2 \tfrac\del{\del u^j}\tfrac\del{\del u^k} 
+ (\tfrac\del{\del u^j} {c^{jk}}\ord2)\tfrac\del{\del u^k}\Bigr]U_0\\
W_6 &=  - \frac12\tfrac\del{\del u^j} [a\ord2 + \tr c\ord2] \tfrac\del{\del u^j} U_0
\end{split}\end{equation*}

By (\ref{u0def}) we compute
$$
W_1 = [-\del_s^2 + 2ia\ord2 \del_s ] U_0 = \Bigl[\frac{n^2}4 - \frac{nz^2}4 \Bigr] U_0
$$
The contribution to $\sigma$ from $W_1$ is thus:
\begin{equation}\label{cont1}
- \frac{n^2}4
\end{equation}

For $W_2$, from the calculations of $a\ord2$ and $a\ord4$ we have 
\begin{equation*}\begin{split}
- a\ord4 + (a\ord2)^2 &= \frac18 (\del_j \del_l \omega_{mk}) z^j z^l z^m (Jz)^k
- \frac1{24} (\del_k F^i_{jl}) (Jz)^k z^j z^l (\omega z)_i \\
&\qquad - \frac1{12} (\del_k F^i_{jl}) z^k z^j (Jz)^l (\omega z)_i
- \frac18 (\del_k F^i_{jl}) z^k z^j z^l z_i \\
&\qquad- \frac18  (\del_j\del_k \beta_{lm}) (Jz)^l z^j z^k (Jz)^m 
\end{split}\end{equation*}
We symmetrize and take the contractions to find the contribution to $\sigma$:
\begin{equation*}\begin{split}
&\frac18 (\del^j \del_j \omega_{mk}) \omega^{mk} + \frac14 (\del_j \del^l \omega_{lk})
\omega^{jk} - \frac1{12} (\beta^{lm} \del_k F^k_{lm}) - \frac1{6} (\del^k F^l_{kl}) \\
&\qquad - \frac14 (\del_j\del_k \beta_{lm}) \omega^{jl} \omega^{km}  -\frac18 (\beta^{lm}\del^k \del_k
\beta_{lm})
\end{split}\end{equation*}

Let us simplify this expression.  By $d\bar \omega =0$ we have
$$
(\del_j \del^l \omega_{lk}) \omega^{jk} = \frac12 (\del^j \del_j \omega_{mk}) \omega^{mk}.
$$
From $\bar\omega_{mk} = -\bar\beta_{mr} \bar J_k^r$ we derive
$$
(\del^j \del_j \omega_{mk}) \omega^{mk} =  \beta^{lm}\del^j\del_j\beta_{lm} 
- (\del^j \del_j J^m_k) J^k_m 
$$
Finally from $\bar J^2=-1$ we obtain
$$
(\del^j \del_j J^m_k) J^k_m = - (\del_j J^m_k) (\del^j J^k_m) = |\nabla J|^2.
$$
Combining these facts gives
$$
\frac18 (\del^j \del_j \omega_{mk}) \omega^{mk} + \frac14 (\del_j \del^l \omega_{lk})
\omega^{jk} = \frac14 \beta^{lm}\del^j\del_j\beta_{lm} - \frac14 |\nabla J|^2
$$

Evaluating the Christoffel symbols gives
\begin{equation*}\begin{split}
\beta^{lm} \del_k F^k_{lm} &= \frac12 \beta^{lm} \del^k [\del_l \beta_{mk} + \del_m \beta_{lk} -
\del_k\beta_{lm}] \\ 
&=   \del^k\del^l \beta_{lk} - \frac12 \beta^{lm}\del^k\del_k\beta_{lm}
\end{split}\end{equation*}
and
$$
\del^k F^l_{kl} = \frac12 \beta^{lm} \del^k\del_k\beta_{lm}
$$
Thus the final contribution from $W_2$ to $\sigma$ is
\begin{equation}\label{cont2}
- \frac14 |\nabla J|^2 - \frac14 (\del_j\del_k \beta_{lm}) \omega^{jl} \omega^{km} 
+\frac1{12} \beta^{lm} \del^k\del_k\beta_{lm} - \frac1{12} \del^j\del^l \beta_{jl}
\end{equation}

By our calculations,
$$
(b\ord2)^2 = \frac19 (\del_l\omega_{km}) z^k z^l (\del_i J_j^m) z^i z^j ,
$$
which (recalling that $\del^jJ_j^m=0$) gives a contribution from $W_3$ of
$$
\frac19 (\del_l \omega_{km}) (\del^k \omega^{lm})
+ \frac19 |\nabla J|^2
$$
By $d\bar\omega = 0$, we have
$$
(\del_l \omega_{km}) (\del^k \omega^{lm}) = - \frac12 (\del_k \omega_{ml}) (\del^k \omega^{lm}) = +\frac12
|\nabla J|^2.
$$
So the contribution from $W_3$ simplifies to 
\begin{equation}\label{cont3}
\frac1{6} |\nabla J|^2
\end{equation}

The terms in $W_4$ are purely imaginary and therefore must contribute zero because $\sigma$ is real.
This can easily be confirmed explicitly.

To compute $W_5$ we need to consider
$$
{c^{jk}}\ord2 \del_{u^j}\del_{u^k}U_0 + (\del_{u^j} {c^{jk}}\ord2)\del_{u^k}U_0
$$
Noting that $\del_{u^j} U_0 = -\frac{u_j}2 U_0$, this becomes
$$
\Bigl[\frac14 c_{jk}\ord2 u^j u^k - \frac12 \beta^{jk} c_{jk}\ord2 - \frac12 u_k (\del_{u^j}
{c^{jk}}\ord2) \Bigr] U_0
$$
If ${c^{jk}}\ord2$ is written $E^{jk}_{lm} u^l u^m$, then under contraction the contribution is
\begin{equation*}\begin{split}
&\frac14 (\beta ^{lm} \beta_{jk} E^{jk}_{lm}+E^{jk}_{jk}+E^{jk}_{kj}) - \frac12 \beta ^{lm} \beta_{jk}
E^{jk}_{lm} - \frac12 (E^{jk}_{jk}+E^{jk}_{kj})\\ 
&\qquad= -\frac14 (\beta ^{lm} \beta_{jk} E^{jk}_{lm}+E^{jk}_{jk}+E^{jk}_{kj})
\end{split}\end{equation*}
This is the same as the contribution of 
$$
-\frac14 c_{jk}\ord2 u^j u^k  = -\frac18 (\del_j\del_k \beta_{lm}) z^j z^k z^l z^m
+ \frac14 (\del_m F_{jlk}) z^m z^k z^j z^l,
$$
yielding
$$
- \frac18 \beta^{lm} (\del^j\del_j \beta_{lm}) - \frac14 (\del^j\del^k \beta_{jk}) + \frac14 \beta^{lm}\del_k
F^k_{lm} + \frac12 (\del^m F^k_{mk}),
$$
which vanishes upon substitution of the $F$.  Hence the total contribution of $W_5$ to $\sigma$ is zero.

Finally, we evaluate the expression appearing in $W_6$:
\begin{equation*}\begin{split}
\frac14 u^j \del_{u^j} [a\ord2 + \tr c\ord2] &= \frac12 [a\ord2 + \tr c\ord2] \\
&=  \frac14 (\beta^{lm}\del_j\del_k\beta_{lm}) z^j z^k
- \frac16 (\del_l F^l_{ik}) z^i z^k - \frac13 (\del_m F^l_{il}) z^m z^i
\end{split}\end{equation*}
The contribution is
$$
\frac14 (\beta^{lm}\del^k\del_k\beta_{lm}) 
- \frac16 (\beta^{ik}\del_l F^l_{ik}) - \frac13 (\del^m F^l_{ml}).
$$
This contribution from $W_6$ can be reduced to
\begin{equation}\label{cont4}
\frac16(\beta^{lm}\del^k\del_k\beta_{lm}) - \frac16 (\del^k\del^l\beta_{kl})
\end{equation}

Adding together (\ref{cont1}), (\ref{cont2}), (\ref{cont3}), and (\ref{cont4}) gives
$$
\sigma = - \frac{n^2}4 - \frac1{12} |\nabla J|^2 - \frac14 (\del_j\del_k \beta_{lm}) \omega^{jl} \omega^{km} 
+\frac1{4} \beta^{lm} \del^k\del_k\beta_{lm} - \frac1{4} \del^j\del^l \beta_{jl}
$$
The last three terms on the right-hand side could be written in terms of the curvature tensors:
$$
- \frac14 (\del_j\del_k \beta_{lm}) \omega^{jl} \omega^{km} 
+\frac1{4} \beta^{lm} \del^k\del_k\beta_{lm} - \frac1{4} \del^j\del^l \beta_{jl}
= \frac14(R + \frac12 R_{ljkm} \omega^{lj} \omega^{km}).
$$
To complete the calculation we cite a lemma which can be found, for example, in \cite{Hs}.
\begin{lemma}
For an almost K\"ahler manifold, 
$$
R + \frac12 R_{ljkm} \omega^{lj} \omega^{km} = -\frac12 |\nabla J|^2.
$$
\end{lemma}

This lemma leads us to the final result that
\begin{equation}
\sigma = - \frac{n^2}4 - \frac5{24} |\nabla J|^2
\end{equation}

\subsection{Quasimodes}
Let us introduce the function 
$$
h(x) = - \frac5{24} |\nabla J(x)|^2
$$

\begin{proposition}\label{psiseq}
Fix $x_0 \in X$ and let $\Gamma=\pi^{-1}(x_0)$. 
There exists a sequence $\psi_k \in L^2(Z)_k$ with $\norm{\psi_k} = 1$ such that
\begin{equation}\label{psiest}
\norm{(\Delta_h - nk - h(x_0)) \psi_k} = O(k^{-1/2}).
\end{equation} 
Moreover, $\psi_k$ is asymptotically localized on $\Gamma$ in the sense that if $\varphi\in \cinf(Z)$
vanishes to order $m$ on $\Gamma$, then
\begin{equation}\label{psiconc}
\pair{\psi_k, \varphi \psi_k} = O(k^{-m/2}).
\end{equation} 
\end{proposition}
\begin{proof}
Let $W$ be a neighborhood of $\Gamma$ in which Fermi coordinates $(s,y)$ are valid, and $\chi \in
\cinf(Z)$ a cutoff function with supp$(\chi)\subset W$ and $\chi=1$ in some neighborhood of $\Gamma$.  
Then we define the sequence $\psi_k \in \cinf(Z)_k$ by
$$
\psi_k(s,y) = \Lambda_k \chi e^{i\kappa s}[U_0 + \kappa^{-1}U_1],
$$
where $U_j(s,y)$ are the solutions obtained above, $\kappa = k + n/2$, and $\Lambda_k$ normalizes
$\norm{\psi_k}=1$.  This could be written as
\begin{equation}\label{pskpoly}
\psi_k(s,y) = \Lambda_k \chi e^{iks}[P_0 + P_2(y) + \kappa P_4(y)]e^{-\kappa y^2/4},
\end{equation}
where $P_l$ is a polynomial of degree $l$ (with coefficients independent of $k$).  Since $P_0 = 1 + O(k^{-1})$,
we have that
$$
\Lambda_k \sim \Bigl(\frac{k}{2\pi}\Bigr)^{n/2} \text{ as }k\to\infty.
$$
The concentration of $\psi_k$ on $\Gamma$ described in (\ref{psiconc}) then follows immediately
from (\ref{pskpoly}).

By virtue of the factor $e^{-\kappa y^2/4}$, we can turn the formal considerations used to obtain the operators
$\calL_j$ into estimates.  With cutoff, $\chi \calL_j$ could be considered an operator on $Z$ with support
in $W$.  By construction we have
$$
\chi \Bigl[ e^{- i\kappa s}\Delta_Z e^{i\kappa s} - \kappa^2 - \kappa\calL_0 - \sqrt\calL_1 - \calL_2  \Bigr]
= \sum_{l,m,|\beta|\le2} E_{l,m,\beta}(s,y) \kappa^l\del_s^m \del_{y}^\beta,
$$
where $A_{l,m,\beta}$ is supported in $W$ and vanishes to order $2l+|\beta|+1$ at $y=0$.
We also have
$$
(\kappa\calL_0 + \sqrt{\kappa}\calL_1 + \calL_2 - \sigma)(U_0 + \kappa^{-1} U_1)
= \kappa^{-1} (\sqrt{\kappa} \calL_1 +\calL_2 - \sigma) U_1 
$$
Combining these facts with the definition of $\psi_k$ we deduce that
$$
(\Delta_Z-\kappa^2-\sigma) \psi_k(s,y) = \Lambda_k \sum_{l\le 4} k^l  F_{l}(s,y) e^{-\kappa y^2/4},
$$
where $F_l$ is supported in $W$ and vanishes to order $2l+1$ at $y=0$.  
Using this order of vanishing we estimate
$$
\norm{\Lambda_k k^l F_{l} e^{-\kappa y^2/4}}^2 = O(k^{-1}).
$$
Noting that $\Delta_Z-\kappa^2-\sigma = \Delta_h - nk - h(x_0)$ on $L^2(Z)_k$, we obtain 
the estimate (\ref{psiest}).
\end{proof}

\section{Spectral density function}

Let $\psi_k\in L^2(Z)_k$ be the sequence produced by Proposition \ref{psiseq}. 
As in \S2, we let $\Pi_k$ denote the orthogonal projection onto the span of 
low-lying eigenvectors of $\Delta_h - nk$.  Consider 
$$
\phi_k = \Pi_k\psi_k \qquad \eta_k = (I-\Pi_k)\psi_k.
$$
By Theorem \ref{GUthm} (for $k$ sufficiently large, which we'll assume throughout),
$$
\norm{(\Delta_h-nk)\phi_k} < M,\qquad  \norm{(\Delta_h-nk)\eta_k} > ak\norm{\eta_k}.
$$ 
By Proposition \ref{psiseq} we have a uniform bound
$$
\norm{(\Delta_h-nk) \psi_k} \le C,
$$ 
so these estimates imply in particular that
$$
ak\norm{\eta_k} < C+M.
$$ 
Hence $\norm{\eta_k} = O(k^{-1})$.

From Lemma \ref{qlemma} we know that $q$ satisfies
$$
\pair{\phi_k, (\Delta_h-nk-\pi^*q)\phi_k} = O(1/k).
$$  
Let $r_k = (\Delta_h - nk + h(x_0))\psi_k$, which by Proposition \ref{psiseq} satisfies
$\norm{r_k} = O(k^{-1/2})$.  So
\begin{equation}\label{phiphi}\begin{split}
&\pair{\phi_k, (\Delta_h-nk-\pi^*q)\phi_k} \\
&\quad= \pair{\phi_k,(h(x_0)-\pi^*q) \phi_k} + \pair{\phi_k,
(\Delta_h-nk-h(x_0))\phi_k}\\
&\quad= \pair{\phi_k,(h(x_0)-\pi^*q) \phi_k} + \pair{\phi_k, r_k}
 - \pair{\phi_k, (\Delta_h-nk-h(x_0))\eta_k}.
\end{split}\end{equation}
The left-hand side is $O(1/k)$, while the second term on the right is $O(k^{-1/2})$, 
The third term term on the right-hand side is equal to
$$
\pair{(\Delta_h-nk)\phi_k,\eta_k} < M\norm{\eta_k} = O(k^{-1}).
$$  
Therefore, the first term on the right-hand side of (\ref{phiphi}) can be estimated
$$
\pair{\phi_k,(h(x_0)-\pi^*q) \phi_k} = O(k^{-1/2}).
$$
Because $\norm{\eta_k} = O(1/k)$ this implies also that 
$$
h(x_0) - \pair{\psi_k, (\pi^*q)\psi_k} = O(k^{-1/2}).
$$
Since $q$ is smooth, the localization of $\psi_k$ on $\Gamma$ from Proposition \ref{psiseq} 
implies that 
$$
\pair{\psi_k, (\pi^*q)\psi_k} = q(x_0) + O(k^{-1/2}).
$$
Thus $q(x_0) = h(x_0)$.  This proves Theorem \ref{mainthm}.

\end{document}